\documentclass[pdflatex,sn-mathphys-num]{sn-jnl}

\usepackage{graphicx}%
\usepackage{multirow}%
\usepackage{amsmath,amssymb,amsfonts}%
\usepackage{amsthm}%
\usepackage{mathrsfs}%
\usepackage[title]{appendix}%
\usepackage{xcolor}%
\usepackage{subcaption}
\usepackage{textcomp}%
\usepackage{manyfoot}%
\usepackage{booktabs}%
\usepackage{algorithm}%
\usepackage{algorithmicx}%
\usepackage{algpseudocode}%
\usepackage{listings}%

\numberwithin{equation}{section}

\usepackage{geometry}
\usepackage[T1]{fontenc}
\usepackage{lmodern}

\usepackage{tikz}
\usepackage{pgfplots}
\pgfplotsset{compat=1.18}

\usetikzlibrary{calc,shadows.blur,positioning,arrows.meta}
\usetikzlibrary{matrix, fit}

\definecolor{PlotBlack}{RGB}{0,0,0}
\definecolor{AxisGray}{RGB}{60,60,60}

\colorlet{BlueC}{blue!80}
\colorlet{RedC}{red!80}
\colorlet{YellowC}{yellow!90!orange}
\colorlet{GreenC}{green!65!black}

\newcommand{\ULHLCellInset}[7]{%
  \path[
    fill=#4, fill opacity=#5,
    draw=#4!70!black, draw opacity=0.95,
    line width=#7
  ]
    ({#2 + #6},{#1 - #3 + 1 + #6}) rectangle ++({1 - 2*#6},{1 - 2*#6});
}

\newcommand{\QAttackSwatch}[1]{%
  {\setlength{\fboxsep}{0pt}%
   \fcolorbox{#1!70!black}{#1}{\rule{0pt}{1.6ex}\hspace{1.6ex}}}%
}

\newcommand\QueensChessBoardpAttack[5][]{%
  \def\N{#2}%
  \def\positions{#3}%
  \def\highlights{#4}%
  \def\LegendOpt{#1}%
  \begin{tikzpicture}[scale=#5, every node/.style={minimum size=1cm * #5, inner sep=0pt}]
    \foreach \i in {1,...,\N}{
      \foreach \j in {1,...,\N}{
        \pgfmathparse{mod(\i+\j,2) ? "brown!80" : "cream!70"}%
        \edef\color{\pgfmathresult}%
        \path[fill=\color] (\i,\j) rectangle ++(1,1);
      }
    }

    \def\op{0.92}%
    \def\lw{1.05pt}%

    \foreach \qx/\qy in \positions{
      \foreach \r in {1,...,\N}{
        \ifnum\r=\qy\relax\else
          \ULHLCellInset{\N}{\qx}{\r}{BlueC}{\op}{0.02}{\lw}%
        \fi
      }
    }

    \foreach \qx/\qy in \positions{
      \foreach \c in {1,...,\N}{
        \ifnum\c=\qx\relax\else
          \ULHLCellInset{\N}{\c}{\qy}{RedC}{\op}{0.08}{\lw}%
        \fi
      }
    }

    \foreach \qx/\qy in \positions{
      \pgfmathtruncatemacro{\d}{\qx-\qy}%
      \foreach \r in {1,...,\N}{
        \ifnum\r=\qy\relax\else
          \pgfmathtruncatemacro{\c}{\r+\d}%
          \ifnum\c>0\relax
            \ifnum\c<\numexpr\N+1\relax
              \ULHLCellInset{\N}{\c}{\r}{YellowC}{\op}{0.14}{\lw}%
            \fi
          \fi
        \fi
      }
    }

    \foreach \qx/\qy in \positions{
      \pgfmathtruncatemacro{\s}{\qx+\qy}%
      \foreach \r in {1,...,\N}{
        \ifnum\r=\qy\relax\else
          \pgfmathtruncatemacro{\c}{\s-\r}%
          \ifnum\c>0\relax
            \ifnum\c<\numexpr\N+1\relax
              \ULHLCellInset{\N}{\c}{\r}{GreenC}{\op}{0.20}{\lw}%
            \fi
          \fi
        \fi
      }
    }

    \foreach \x/\y in \highlights{
      \ULHLCellInset{\N}{\x}{\y}{black}{0.35}{0.26}{0.9pt}%
    }

    \def\temp{}%
    \ifx\LegendOpt\temp\relax
    \else
      \node[
        anchor=north west,
        draw=black!25,
        fill=white,
        rounded corners=2pt,
        inner sep=5pt,
        align=left,
        font=\scriptsize
      ] at (\N+1.2,\N+1.05) {%
        \textbf{Attacks}\\[2pt]
        \begin{tabular}{@{}l l@{}}
          \QAttackSwatch{BlueC}   & Column\\
          \QAttackSwatch{RedC}    & Row\\
          \QAttackSwatch{YellowC} & Diagonal\\
          \QAttackSwatch{GreenC}  & Anti-diagonal\\
        \end{tabular}%
      };
    \fi

    \foreach \x/\y in \positions{
      \node at (\x+.5, \N-\y+1.5)
        {\includegraphics[width=0.35cm]{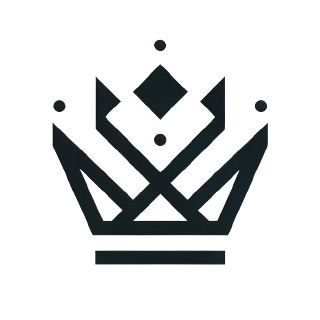}};
    }
  \end{tikzpicture}%
}

\newcommand{\CapSwatch}[1]{%
  {\setlength{\fboxsep}{0pt}%
   \raisebox{0.15ex}{\fcolorbox{#1!70!black}{#1}{\rule{0pt}{1.35ex}\hspace{1.35ex}}}}%
}

\definecolor{cream}{RGB}{255,253,208}
\theoremstyle{thmstyleone}%
\newtheorem{theorem}{Theorem}
\newtheorem{proposition}[theorem]{Proposition}%

\theoremstyle{thmstyletwo}%

\theoremstyle{thmstylethree}%

\begin{document}

\title[Improved asymptotic upper bound on the $n$-queens completion threshold]{Improved asymptotic upper bound on the $n$-queens completion threshold}


\author*[1]{\fnm{Hugo Møller} \sur{Nielsen}}\email{hugomn2002@gmail.com}

\affil[1]{\orgdiv{Department of Applied Mathematics and Computer Science}, \orgname{Technical University of Denmark (DTU)},
	\orgaddress{\city{Kgs. Lyngby}, \postcode{2800}, \country{Denmark}}}


\abstract{
	The $n$-queens completion threshold $qc(n)$ is the largest integer $k < n$ such that any
	placement of $k$ mutually non-attacking queens on an $n \times n$ chessboard
	can be completed to an $n$-queens configuration by adding $n - k$ queens.
	For all sufficiently large $n$, we improve the previously best-known upper bound on
	$qc(n)$ from $qc(n) \leq 0.241n$ to $qc(n) \leq 0.216n$,
	by constructing a non-completable partial configuration of fewer than $0.216n$ queens.
}




\maketitle

\section{Introduction}
An $n$-queens configuration is a placement of $n$ mutually non-attacking queens on an $n \times n$ chessboard.
That is, no two queens share a row, column, or diagonal. The $n$-queens problem,
posed by Max Bezzel in 1848, asks: for a given board size $n$, how many distinct $n$-queens configurations
exist? We denote this number by $Q(n)$.

An interesting variant of the $n$-queens problem, posed by Franz Nauck in 1850, is the $n$-queens completion problem, which asks
whether a partial $n$-queens configuration, i.e., a set of $k \leq n$ mutually non-attacking queens,
is completable in the sense that $n - k$ queens can be added to make it a valid $n$-queens
configuration \cite{Campbell1977,Gent2017NQueensCompletion}. This problem is NP-complete \cite{Gent2017NQueensCompletion},
but in some cases the completability of a partial configuration can still be determined easily.
For instance, since $n$-queens configurations exist for all $n \geq 4$ \cite{Pauls1874Damenproblem},
the empty partial configuration is completable whenever $n \geq 4$.
For $n \in \mathbb{N}$, the $n$-queens completion threshold $qc(n)$ is the largest integer $k < n$ such
that every partial $n$-queens configuration consisting of $k$ or fewer queens can be completed to a
full $n$-queens configuration, by adding $n - k$ additional queens to the partial configuration.
In recent work, Glock et al.\ \cite{Glock2022NQueens} showed that asymptotically $qc(n) \geq n/60$
by relating completion to rainbow matchings in an associated bipartite graph, and proved the upper bound
$qc(n) \leq 0.241n$ by giving an explicit non-completable construction.
In this work, we improve this upper bound by showing that the construction in \cite{Glock2022NQueens} is
non-completable once $0.216n$ queens have been placed, and thus $qc(n) \leq 0.216n$.

\section{Background}
We briefly outline the notation used in this work and some basic properties of the n-queens completion problem,
and use a similar setup as in \cite{Glock2022NQueens}.

\subsection{Representing $n$-queens configurations}
Let $n$ denote the board size throughout.
We represent the chessboard as the two-dimensional grid $[n] \times [n]$,
where $[n] = \lbrace 1, 2, \dots, n \rbrace$.
We write the rows, columns, diagonals and anti-diagonals as
\begin{align*}
	R_{i}     & = \lbrace (i,j) | j \in [n] \rbrace,                        \\
	C_{j}     & = \lbrace (i,j) | i \in [n] \rbrace,                        \\
	D_{k}^{-} & = \lbrace (i,j) \in [n] \times [n] | j - i + n = k \rbrace, \\
	D_{k}^{+} & = \lbrace (i,j) \in [n] \times [n] | i + j - 1 = k \rbrace.
\end{align*}
The set of all $n$ rows, $n$ columns, $2n-1$ diagonals and $2n-1$ anti-diagonals on an $n \times n$ chessboard are
denoted by respectively $\mathcal{R}_n, \mathcal{C}_n, \mathcal{D}^{-}_{n}$ and $\mathcal{D}^{+}_{n}$.
The set of all lines on an $n \times n$ chessboard is denoted by
$\mathcal{L}_n = \mathcal{R}_n \cup \mathcal{C}_n \cup \mathcal{D}^{+}_{n} \cup \mathcal{D}^{-}_{n}$.
We identify the square in row $i \in [n]$ and column $j \in [n]$ with the coordinate $(i,j) \in [n] \times [n]$,
and define an $n$-queens configuration to be a set $Q \subseteq [n] \times [n]$ of cardinality $|Q| = n$ such that
\begin{align*}
	 & |Q \cap R_i| \leq 1 \;\;\; \forall i \in [n],          &  & |Q \cap C_j| \leq 1 \;\;\; \forall j \in [n],          \\ 
	 & |Q \cap D_{k}^{-}| \leq 1 \;\;\; \forall k \in [2n-1], &  & |Q \cap D_{k}^{+}| \leq 1 \;\;\; \forall k \in [2n-1]. 
\end{align*}
Given a set $S \subseteq [n] \times [n]$, we call $U \subseteq [n] \times [n]$ a completion of $S$ if
$S \cup U$ is an $n$-queens configuration and $S \cap U = \emptyset$.

\subsubsection{Weighting \& Covering}
Let $w : \mathcal L_n\to[0,1]$ be a weighting function. For a set of lines $\mathcal{L} \subseteq \mathcal{L}_n$,
define its total weight by
\begin{align*}
	w(\mathcal L) = \sum_{L\in\mathcal L} w(L).
\end{align*}
We say that $w$ \emph{covers} a square $(i,j)\in[n]\times[n]$ if the total weight of the four lines through $(i,j)$
is at least $1$, that is,
\begin{align*}
	w(R_i)+w(C_j)+w(D_{i+j-1}^{+})+w(D_{j-i+n}^{-})\ge 1.
\end{align*}
A set $\Lambda \subseteq [n] \times [n]$ is covered if every square in $\Lambda$ is covered.
We will need the following result from \cite{Glock2022NQueens}:
\begin{proposition} \label{prop:weighting_method}
	Let $Q^\prime$ be a partial $n$-queens configuration and let $\Lambda \subseteq [n] \times [n]$ be
	the set of squares not attacked by $Q^\prime$. If there exists a line weighting which covers $\Lambda$
	and has $w(\mathcal{L}_n) < n - |Q^\prime|$, then $Q^\prime$ cannot be completed.
\end{proposition}
For any specific partial $n$-queens configuration $Q^\prime$, Proposition~\ref{prop:weighting_method} naturally gives rise to the problem of finding
a minimal covering weighting function for $\Lambda \subseteq [n] \times [n]$, the set of squares not attacked by any queen in $Q^\prime$.

\subsubsection{Finding the optimal weighting function $w : \mathcal{L}_n \to [0,1]$} \label{sec:optimal_weighting}
Given a partial $n$-queens configuration $Q^\prime \subseteq [n] \times [n]$, we can formulate the problem of finding a minimum-weight cover as a linear program.
The objective is to minimize the total weight $w(\mathcal{L}_n)$ induced by the weighting function, which can be written as
\begin{align} \label{eq:sum_over_lines}
	\min \sum_{i=1}^{n} \left(w(R_i) + w(C_i)\right) + \sum_{k=1}^{2n-1}\left(w(D_{k}^{-}) + w(D_{k}^{+})\right).
\end{align}
Furthermore, we want all non-threatened squares $\Lambda \subseteq [n] \times [n]$ to satisfy the covering constraint, which we can express as such:
\begin{align}
	w(R_i) + w(C_j) + w(D_{j-i+n}^{-}) + w(D_{i+j-1}^{+}) \geq 1 \;\; \forall (i,j) \in \Lambda. \label{eq:covering_constraint}
\end{align}
This yields a linear program with $6n-2$ variables and $|\Lambda| = O(n^2)$ constraints.
Since linear programs can be solved in polynomial time, we can compute a minimum-weight cover for $\Lambda$
even for moderately large $n$.

\section{Improved upper bound}
In this section, we show that $qc(n) \leq 0.216n$ for $n$ large enough. To this end, fix odd $n \in \mathbb{N}$.

\subsection{Non-completable partial configuration}
\begin{figure}[h!]
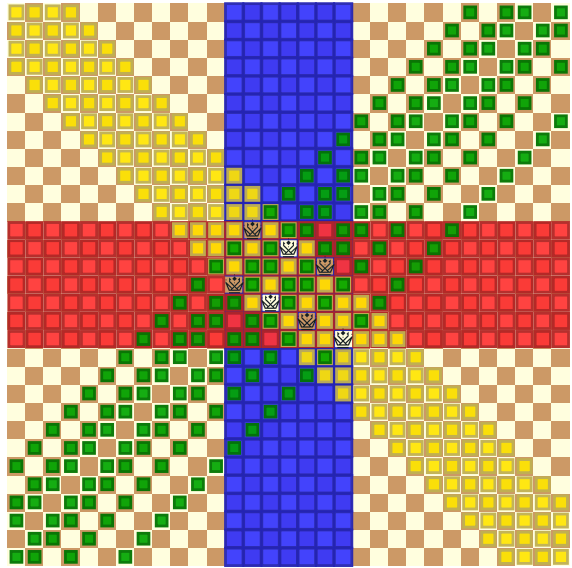

	\centering
	\QueensChessBoardpAttack{31}{14/13, 16/14, 18/15, 13/16, 15/17, 17/18, 19/19}{}{0.24}
	\caption{Partial configuration $Q^\prime$ for $n=31$ and $m=7$.
		Attacked squares are colored by line type:
		\CapSwatch{BlueC} columns,
		\CapSwatch{RedC} rows,
		\CapSwatch{YellowC} diagonals,
		\CapSwatch{GreenC} anti-diagonals.
	}
	\label{fig:partial_n31_m7}
\end{figure}

Let $m \in \mathbb{N}$ with $m \leq n$ and $m \equiv 1 \pmod{6}$ and construct the partial $n$-queens configuration
\begin{align}
	Q^\prime & = \left\lbrace \left(i + \frac{n - m}{2}, 2i + \frac{n-m}{2} \right) \middle| i \in \left\lbrace 1, 2, \dots, \left\lfloor\frac{m}{2}\right\rfloor \right\rbrace \right\rbrace                                                                  \\
	         & \cup \left\lbrace \left(i + \frac{n - m}{2}, 2\left(i - \left\lfloor\frac{m}{2}\right\rfloor\right) - 1 + \frac{n-m}{2} \right) \middle| i \in \left\lbrace \left\lfloor\frac{m}{2}\right\rfloor + 1, \left\lfloor\frac{m}{2}\right\rfloor + 2,
	\dots, m \right\rbrace \right\rbrace. \label{eq:partial_grock}
\end{align}
This is an embedding of a standard $m$-queens construction into the center of an $n \times n$ chessboard
\cite{bell2009constructing,Glock2022NQueens}, see Figure~\ref{fig:partial_n31_m7}.

\begin{figure}[ht!]
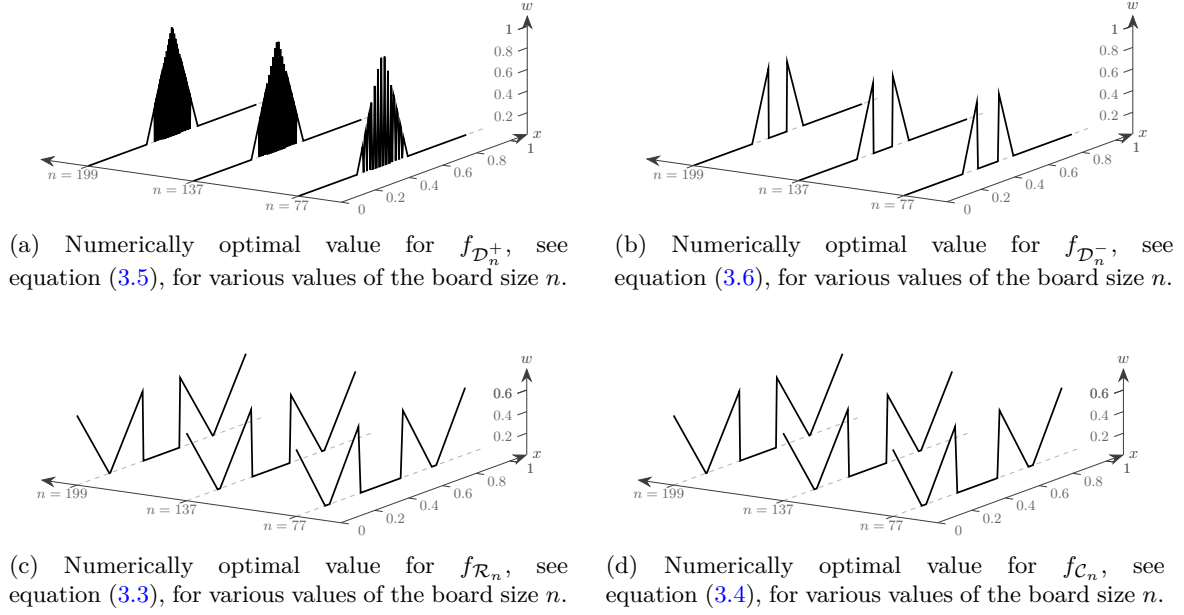

	\centering
	\begin{subcaptiongroup}
		\begin{subfigure}{0.48\textwidth}
			\centering

			\caption{Numerically optimal value for $f_{\mathcal{D}^{+}_{n}}$, see equation~(\ref{eq:fDp_function}), for various values of the board size $n$.}
			\label{fig:Dp_diag_numerical}
		\end{subfigure}
		\hfill
		\begin{subfigure}{0.48\textwidth}
			\centering
			%
			\caption{Numerically optimal value for $f_{\mathcal{D}^{-}_{n}}$, see equation~(\ref{eq:fDm_function}), for various values of the board size $n$.}
			\label{fig:Dm_diag_numerical}
		\end{subfigure}

		\vspace{0.5cm} 

		\begin{subfigure}{0.48\textwidth}
			\centering
			%
			\caption{Numerically optimal value for $f_{\mathcal{R}_{n}}$, see equation~(\ref{eq:fR_function}), for various values of the board size $n$.}
			\label{fig:R_diag_numerical}
		\end{subfigure}
		\hfill
		\begin{subfigure}{0.48\textwidth}
			\centering
			%
			\caption{Numerically optimal value for $f_{\mathcal{C}_{n}}$, see equation~(\ref{eq:fC_function}), for various values of the board size $n$.}
			\label{fig:C_diag_numerical}
		\end{subfigure}
	\end{subcaptiongroup}
	\caption{Optimal weighting function $w$ for various board sizes $n$, the code is available at \cite{nielsen2025qclp}.} \label{fig:optimal_w_numerically}
\end{figure}

We search over admissible values of $m$ and, for each $m$, compute a minimum-weight cover as described in Section~\ref{sec:optimal_weighting}.
Let $t=m/n$. We record the smallest value of $t \in [0,1]$ for which Proposition~\ref{prop:weighting_method} guarantees that $Q^\prime$ is non-completable,
and we store the corresponding optimal weighting $w : \mathcal{L}_n \to [0,1]$.
Doing this for various $n$ and plotting the functions
\begin{align}
	f_{\mathcal{R}_{n}}     & : [n] \to [0,1], \; f_{\mathcal{R}_{n}}(i) = w(R_{i}), \label{eq:fR_function}             \\
	f_{\mathcal{C}_{n}}     & : [n] \to [0,1], \; f_{\mathcal{C}_{n}}(i) = w(C_{i}), \label{eq:fC_function}             \\
	f_{\mathcal{D}_{n}^{+}} & : [2n-1] \to [0,1], \; f_{\mathcal{D}_{n}^{+}}(i) = w(D_{i}^{+}), \label{eq:fDp_function} \\
	f_{\mathcal{D}_{n}^{-}} & : [2n-1] \to [0,1], \; f_{\mathcal{D}_{n}^{-}}(i) = w(D_{i}^{-}), \label{eq:fDm_function}
\end{align}
where the points $(i,f(i))$ are joined by straight line segments, and the $x$-axis rescaled to $[0,1]$,
we get the plots in Figure~\ref{fig:optimal_w_numerically}.

Figure~\ref{fig:optimal_w_numerically} suggests that the optimal weighting $w : \mathcal{L}_n \to [0, 1]$,
viewed via the four normalized functions in (\ref{eq:fR_function}), (\ref{eq:fC_function}), (\ref{eq:fDm_function})
and (\ref{eq:fDp_function}), is approximately independent of $n$ after linear interpolation and rescaling of
the domain.
Motivated by Figure~\ref{fig:optimal_w_numerically}, we propose an explicit candidate for an optimal
weighting function for the partial configuration in (\ref{eq:partial_grock}). To this end, define
\begin{align} \label{eq:anti_diagonal_Q_prime_identifiers}
	A = \left\lbrace 3\ell + n - m - 1 | \ell \in \left\lbrace 1, 2, \dots, \left\lfloor \frac{m}{2} \right\rfloor \right\rbrace \right\rbrace \cup \left\lbrace 3\ell + n - 2m - 1 | \ell \in \left\lbrace \left\lfloor \frac{m}{2} \right\rfloor + 1, \dots, m \right\rbrace \right\rbrace,
\end{align}
the set of indices $k \in [2n-1]$ such that the anti-diagonal $D_k^+$ is threatened by a queen in the
partial configuration $Q^\prime$ in (\ref{eq:partial_grock}).
For the partial $n$-queen configuration $Q^\prime$, define the weighting function $w: \mathcal{L}_{n} \to [0,1]$ by:
\begin{align}
	w(R_i)           & = \begin{cases}
		                     -\alpha \frac{i-1}{n-1} + \alpha\frac{1 - t}{4}         & 0 \leq \frac{i-1}{n-1} < \frac{1-t}{4},                     \\
		                     \alpha \frac{i-1}{n-1} - \alpha \frac{1-t}{4}           & \frac{1-t}{4} \leq \frac{i-1}{n-1} < \frac{1-t}{2},         \\
		                     0                                                       & \frac{1-t}{2} \leq \frac{i-1}{n-1} < 1 - \frac{1-t}{2},     \\
		                     -\alpha \frac{i-1}{n-1} + \alpha - \alpha \frac{1-t}{4} & 1 - \frac{1-t}{2} \leq \frac{i-1}{n-1} < 1 - \frac{1-t}{4}, \\
		                     \alpha \frac{i-1}{n-1} + \alpha \frac{1-t}{4} - \alpha  & 1 - \frac{1-t}{4} \leq \frac{i-1}{n-1} \leq 1,
	                     \end{cases} \label{eq:wRi} \\
	w(C_j)           & = w(R_j), \label{eq:wCj}                                                                                                \\
	w(D_{j-i+n}^{-}) & = \begin{cases}
		                     0                          & \frac{1}{\alpha} \leq \frac{i-j}{n-1} < 1,             \\
		                     1 - \alpha \frac{i-j}{n-1} & \frac{t}{2} \leq \frac{i-j}{n-1} < \frac{1}{\alpha},   \\
		                     0                          & -\frac{t}{2} \leq \frac{i-j}{n-1} < \frac{t}{2},       \\
		                     1 + \alpha \frac{i-j}{n-1} & -\frac{1}{\alpha} \leq \frac{i-j}{n-1} < -\frac{t}{2}, \\
		                     0                          & -1 \leq \frac{i-j}{n-1} \leq -\frac{1}{\alpha},
	                     \end{cases} \label{eq:wDm}                                   \\
	w(D_{i+j-1}^{+}) & = \begin{cases}
		                     0                                     & 0 \leq \frac{i+j-2}{n-1} < 1 - \frac{1}{\alpha},                      \\
		                     1 + \alpha \frac{i+j-2}{n-1} - \alpha & 1 - \frac{1}{\alpha} < \frac{i+j-2}{n-1} \leq 1 \land i+j-1 \notin A, \\
		                     0                                     & i+j-1 \in A,                                                          \\
		                     1 - \alpha \frac{i+j-2}{n-1} + \alpha & 1 \leq \frac{i+j-2}{n-1} < 1 + \frac{1}{\alpha} \land i+j-1 \notin A, \\
		                     0                                     & 1 + \frac{1}{\alpha} \leq \frac{i+j-2}{n-1} \leq 2,
	                     \end{cases} \label{eq:wDp}
\end{align}
for a constant $\alpha \in \mathbb{R}^{+}$ with $\frac{1}{\alpha} > t$. In Section~\ref{sec:covering_property}
we prove that this weighting covers every square not attacked by any queen in $Q^\prime$.
Moreover, in Section~\ref{sec:covering_property} we show that the weighting is such that it allows us
to improve the best known upper bound of $qc(n)$.

\begin{figure}[ht!]
	\centering
	\begin{subfigure}[t]{0.47\textwidth}
		\centering
		\includegraphics[width=\linewidth]{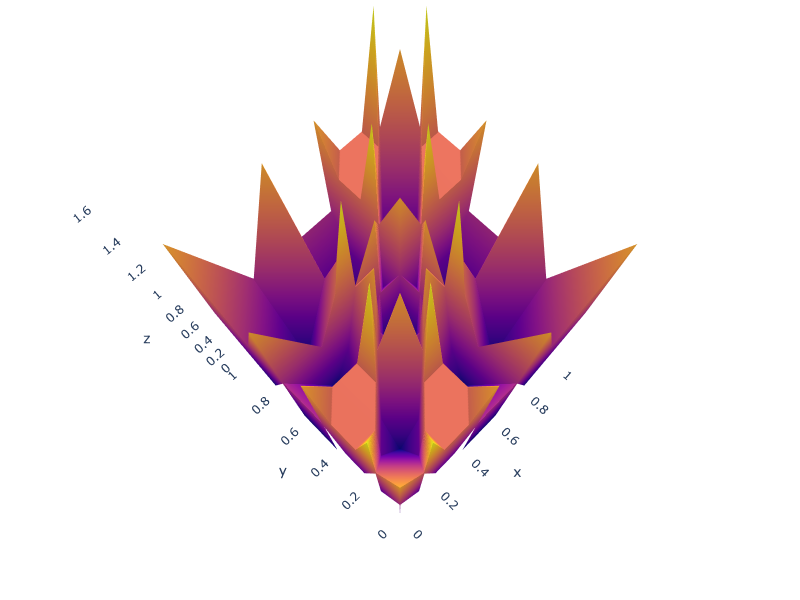}
		\caption{3D plot of the weighting function $w$ for the partial configuration in (\ref{eq:partial_grock}),
			excluding the contribution from $\mathcal{D}_n^+$.
			The values $w(D_k^+)$ are irregular due to the inclusion/exclusion condition involving $A$; see
			Figure~\ref{fig:Dp_diag_numerical}. Here we use $n = 600, t = 0.216$ and $\alpha = 3.29$.}
		\label{fig:3D_plot_RCD}
	\end{subfigure}%
	\hspace{0.04\textwidth} 
	\begin{subfigure}[t]{0.47\textwidth}
		\centering
		\includegraphics[width=0.925\linewidth]{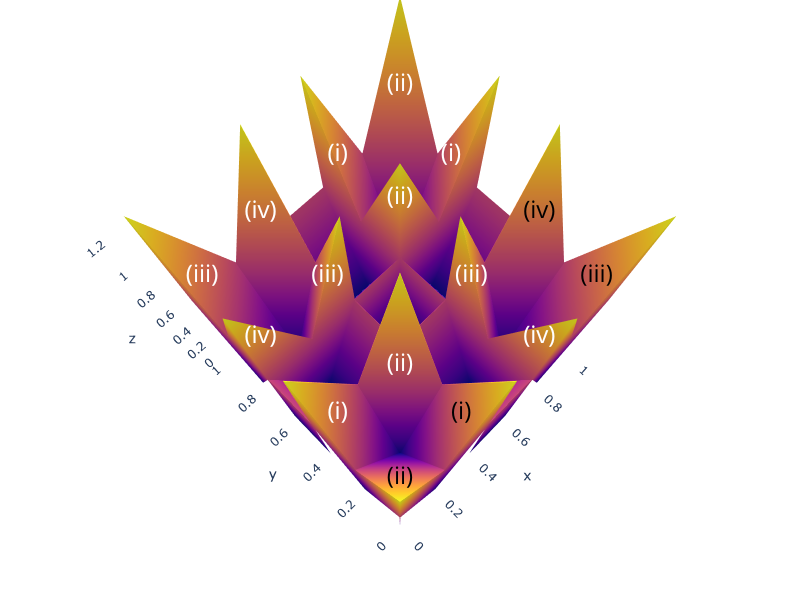}
		\caption{3D plot of the weighting function $w$ for the partial configuration in (\ref{eq:partial_grock}),
			showing only the contribution $w(R_i) + w(C_j)$ for each square $(i,j) \in [n] \times [n]$.
			The black labels indicate the cases considered in Section~\ref{sec:covering_property};
			the corresponding white labels mark the symmetric cases implied by the
			definition of $w : \mathcal{L}_n \to [0,1]$ given in (\ref{eq:wRi}), (\ref{eq:wCj}), (\ref{eq:wDm}),
			and (\ref{eq:wDp}).}
		\label{fig:3D_plot_RC}
	\end{subfigure}
	\caption{
		Density plots of $w(R_i) + w(C_j) + w(D_{j-i+n}^{-})$ and $w(R_i) + w(C_j)$, illustrating the symmetries
		of the weighting function $w : \mathcal{L}_n \to [0,1]$ with respect to the covering constraint (\ref{eq:covering_constraint}).}
	\label{fig:both_figures_3D}
\end{figure}

We visualize $w : \mathcal{L}_{n} \to [0,1]$ by making a density
plot, where square $(i,j)$ is given value $w(R_i) + w(C_j) + w(D^{-}_{j-i+n}) + w(D^{+}_{i+j-1})$ for some value of $n$ and $(i,j) \in [n] \times [n]$,
normalizing the input domain to $[0,1] \times [0,1]$, see Figure~\ref{fig:both_figures_3D}.

Since the functions (\ref{eq:fR_function}), (\ref{eq:fC_function}), and (\ref{eq:fDm_function}) are piecewise linear,
the surfaces in Figure~\ref{fig:3D_plot_RCD} are piecewise planar.
Moreover, notice that $f_{\mathcal{D}^{+}_{n}}$ and $f_{\mathcal{D}^{-}_{n}}$ are defined exactly so as to be zero on the threatened lines and to be precisely
large enough so as to increase the value of $w(R_i) + w(C_j) + w(D_{j-i+n}^{-}) + w(D_{i+j-1}^{+})$ to $1$ for all $(i,j) \in [n] \times [n]$ where
$w(R_i) + w(C_j) < 1$.
In the following sections, we prove that this weighting function covers the squares not attacked by $Q^\prime$, and we find the optimal
value of $\alpha$ that minimizes the total weighting.

\subsection{Covering Property} \label{sec:covering_property}
Fix $m \equiv 1 \pmod{6}$. We show that the weighting $w : \mathcal{L}_n \to [0,1]$ defined in (\ref{eq:wRi}), (\ref{eq:wCj}), (\ref{eq:wDm}), and (\ref{eq:wDp})
covers every square not attacked by any queen in $Q^\prime$.
The squares $(i,j) \in [n] \times [n]$ satisfying $|\frac{i-j}{n-1}| \leq \frac{t}{2}$ or
$|1 - \frac{i+j-2}{n-1}| \leq \frac{1}{\alpha} \land i+j-1 \in A$ are exactly those squares attacked along
diagonals or anti-diagonals by queens in $Q^\prime$.
By the symmetric definition of $w$, see Figure~\ref{fig:optimal_w_numerically} and Figure~\ref{fig:both_figures_3D},
it suffices to consider only the squares $(i,j) \in [n] \times [n]$ satisfying the four cases below.

\paragraph{Case (i): $0 \leq \frac{j-1}{n-1} \leq \frac{1-t}{4} \leq \frac{i-1}{n-1} \leq \frac{1-t}{2}$.}
We split the argument according to the size of $\frac{i - j}{n - 1}$.

First suppose that $\frac{i-j}{n-1} \leq \frac{1}{\alpha}$, then square $(i,j)$ is covered by the inequality:
\begin{align*}
	 & w(R_i) + w(C_j) + w(D_{j-i+n}^{-}) + w(D_{i+j-1}^{+}) \geq w(R_i) + w(C_j) + w(D_{j-i+n}^{-})                                      \\
	 & \geq \alpha\frac{i-1}{n-1} - \alpha \frac{1-t}{4} - \alpha \frac{j-1}{n-1} + \alpha \frac{1-t}{4} + 1 - \alpha\frac{i-j}{n-1} = 1.
\end{align*}
It remains to consider the complementary case $\frac{1}{\alpha} < \frac{i-j}{n-1}$. In this case, square $(i,j)$
is covered by:
\begin{align*}
	 & w(R_i) + w(C_j) + w(D_{j-i+n}^{-}) + w(D_{i+j-1}^{+}) \geq w(R_i) + w(C_j)                                                 \\
	 & = \alpha\frac{i-1}{n-1} - \alpha \frac{1-t}{4} - \alpha\frac{j-1}{n-1} + \alpha \frac{1-t}{4} = \alpha\frac{i-j}{n-1} > 1.
\end{align*}
This proves the covering property in Case (i) and, by symmetry, in the corresponding cases, see Figure~\ref{fig:3D_plot_RC}.

\paragraph{Case (ii): $0 \leq \frac{j-1}{n-1} \leq \frac{i-1}{n-1} \leq \frac{1-t}{4}$.}
For this case we either have
$\frac{i-j}{n-1} \leq \frac{1}{\alpha}$, in which case we cover square $(i,j)$ by the inequality:
\begin{align*}
	 & w(R_i) + w(C_j) + w(D_{j-i+n}^{-}) + w(D_{i+j-1}^{+}) \geq w(R_i) + w(C_j) + w(D_{j-i+n}^{-})                                  \\
	 & \geq - \alpha\frac{i-1}{n-1} + \alpha \frac{1-t}{4} - \alpha\frac{j-1}{n-1} + \alpha \frac{1-t}{4} + 1 - \alpha\frac{i-j}{n-1} \\
	 & = -\alpha \left(\frac{2i-2}{n-1} - \frac{1-t}{2}\right) + 1 \geq 1,
\end{align*}
or we have $\frac{i-j}{n-1} > \frac{1}{\alpha}$ which by the calculations in Case~(i) will cover square $(i,j)$.
This proves the covering property for Case (ii) and, by symmetry, for the corresponding cases.

\paragraph{Case (iii): $1 - \frac{1-t}{4} \leq \frac{i-1}{n-1} \leq 1$ and $0 \leq \frac{j-1}{n-1} \leq \frac{1-t}{4}$ and $i+j-1 \notin A$.}
We split the argument according to the size of $\frac{i+j-2}{n-1}$.

First suppose that $1 - \frac{1}{\alpha} \leq \frac{i+j-2}{n-1} \leq 1$, then the square $(i,j)$ is covered
by the inequality:
\begin{align*}
	 & w(R_i) + w(C_j) + w(D_{j-i+n}^{-}) + w(D_{i+j-1}^{+}) \geq w(R_i) + w(C_j) + w(D_{i+j-1}^{+})                                                   \\
	 & = \alpha\frac{i-1}{n-1} + \alpha \frac{1-t}{4} - \alpha - \alpha \frac{j-1}{n-1} + \alpha \frac{1-t}{4} + 1 + \alpha \frac{i+j-2}{n-1} - \alpha \\
	 & = 1 + \alpha \left(2\frac{i-1}{n-1} + \frac{1-t}{2} - 2\right) \geq 1.
\end{align*}
It remains to consider the complementary case $\frac{i+j-2}{n-1} < 1 - \frac{1}{\alpha}$.
In this case, square $(i,j)$ is covered by:
\begin{align*}
	 & w(R_i) + w(C_j) + w(D_{j-i+n}^{-}) + w(D_{i+j-1}^{+}) \geq w(R_i) + w(C_j) + w(D_{i+j-1}^{+})                                                       \\
	 & = \alpha\frac{i-1}{n-1} + \alpha \frac{1-t}{4} - \alpha - \alpha \frac{j-1}{n-1} + \alpha \frac{1-t}{4}                                             \\
	 & = \alpha \left(-1 + \frac{1-t}{2} + 2\frac{i-1}{n-1} - \frac{j+i-2}{n-1} \right) \geq \alpha \left(-2 + 2\frac{i-1}{n-1} + \frac{1-t}{2}\right) + 1 \\
	 & \geq 1.
\end{align*}
Thus, we have proved the covering property in this case and all symmetric cases.

\paragraph{Case (iv): $1 - \frac{1-t}{4} \leq \frac{i-1}{n-1} \leq 1$ and $\frac{1-t}{4} \leq \frac{j-1}{n-1} \leq \frac{1-t}{2}$ and $i+j-1 \notin A$.}
For this final case we again split the argument according to the size of $\frac{i+j-2}{n-1}$.

First suppose that $1 \leq \frac{i+j-2}{n-1} \leq 1 + \frac{1}{\alpha}$, then square $(i,j)$ is covered by the inequality:
\begin{align*}
	 & w(R_i) + w(C_j) + w(D_{j-i+n}^{-}) + w(D_{i+j-1}^{+}) \geq w(R_i) + w(C_j) + w(D_{i+j-1}^{+})                                                         \\
	 & = \alpha \frac{i-1}{n-1} + \alpha \frac{1-t}{4} - \alpha + \alpha \frac{j-1}{n-1} - \alpha \frac{1-t}{4} + 1 - \alpha \frac{i+j-2}{n-1} + \alpha = 1.
\end{align*}
It remains to consider the complementary case $\frac{i+j-2}{n-1} > 1 + \frac{1}{\alpha}$.
In this case, square $(i,j)$ is covered by:
\begin{align*}
	 & w(R_i) + w(C_j) + w(D_{j-i+n}^{-}) + w(D_{i+j-1}^{+}) \geq w(R_i) + w(C_j)                                                                                 \\
	 & =\alpha \frac{i-1}{n-1} + \alpha \frac{1-t}{4} - \alpha + \alpha \frac{j-1}{n-1} - \alpha \frac{1-t}{4} = \alpha\left(\frac{j+i-2}{n-1} - 1\right) \geq 1.
\end{align*}
Thus, we have proved the covering property in this case and all symmetric cases.

As illustrated in Figure~\ref{fig:3D_plot_RC}, the cases $(i), (ii), (iii),$ and $(iv)$ and
their symmetric counterparts imply that $w$ covers every square not attacked by any queen in $Q^\prime$
on the $n \times n$ board.

\subsection{Total weight of the cover}
As per Proposition~\ref{prop:weighting_method}, we now need to bound the sum in (\ref{eq:sum_over_lines}) for $w$, so as to figure out for
which $t$, equivalently $m$, we are ensured to have a non-completable configuration $Q^\prime$. We bound (\ref{eq:sum_over_lines}) on a
per line-type manner. 
We estimate the sum over the rows and columns using Riemann integration and use Big-$O$ notation with respect to $n$.
As per (\ref{eq:wRi}) and (\ref{eq:wCj}) we obtain   
\begin{align*} 
	\sum_{j = 1}^{n} w(C_j) = \sum_{i = 1}^{n} w(R_i) = 4(n-1)\int_{0}^{\frac{1-t}{4}}-\alpha x + \alpha\frac{1 - t}{4} dx + O(1) = \alpha \frac{(1 - t)^2}{8}(n-1) + O(1).
\end{align*}
The sum over the diagonals is similarly estimated via Riemann-integration by 
\begin{align*}
	\sum_{k = 1}^{2n-1}w(D_k^-) = 2(n-1)\int_{\frac{t}{2}}^{\frac{1}{\alpha}}1 - \alpha x dx + O(1)
	= (n-1) \left(\frac{1}{\alpha} - t + \alpha\frac{t^2}{4}\right) + O(1).
\end{align*}
Finally, to compute the sum over the anti-diagonals we remind the reader of the anti-diagonal identifiers $A$
for elements of $\mathcal{D}_n^{+}$ with a non-empty intersection with $Q^\prime$, see (\ref{eq:anti_diagonal_Q_prime_identifiers}).
We introduce the helper function $\hat{w} : \mathcal{D}_n^+ \to [0,1]$ defined by
\begin{align*}
	\hat{w}(D_{i+j-1}^{+}) = \max\left(0, 1 - \alpha \left|\frac{i+j-2}{n-1} - 1\right|\right) = \begin{cases}
		                                                                                             0                                     & 0 \leq \frac{i+j-2}{n-1} \leq 1 - \frac{1}{\alpha}, \\
		                                                                                             1 + \alpha \frac{i+j-2}{n-1} - \alpha & 1 - \frac{1}{\alpha} \leq \frac{i+j-2}{n-1} \leq 1, \\
		                                                                                             1 - \alpha \frac{i+j-2}{n-1} + \alpha & 1 \leq \frac{i+j-2}{n-1} \leq 1 + \frac{1}{\alpha}, \\
		                                                                                             0                                     & 1 + \frac{1}{\alpha} \leq \frac{i+j-2}{n-1} \leq 2,
	                                                                                             \end{cases} 
\end{align*}
for all $(i,j) \in [n] \times [n]$. By construction,
\begin{align*}
	\sum_{k=1}^{2n-1}w(D_k^{+}) & = \sum_{k=1}^{2n-1} \hat{w}(D_k^{+}) - \sum_{k \in A}\hat{w}(D_k^{+})
	= O(1) + 2(n-1)\int_{1-\frac{1}{\alpha}}^{1}\left(1 + \alpha x - \alpha\right) dx                                                                        \\
	                            & - \sum_{\ell = 1}^{\left\lfloor \frac{m}{2} \right\rfloor}\hat{w}(D_{3\ell + n - m - 1}^+)
	- \sum_{\ell = \left\lfloor \frac{m}{2} \right\rfloor + 1}^{m}\hat{w}(D_{3\ell + n - 2m - 1}^+)
	= O(1) + (n-1) \frac{1}{\alpha}                                                                                                                          \\
	                            & - \sum_{\ell = 1}^{\left\lfloor \frac{nt}{2} \right\rfloor} \left(1 - \alpha \left|\frac{3\ell - m - 1}{n-1}\right|\right)
	- \sum_{\ell = \left\lfloor \frac{nt}{2} \right\rfloor + 1}^{nt} \left(1 - \alpha\left|\frac{3\ell - 2m - 1}{n - 1}\right|\right)                        \\
	                            & = (n-1)\left(\frac{1}{\alpha} - t\right) + \alpha\int_{0}^{\frac{nt}{2}}\left|\frac{3x - nt - 1}{n-1}\right|dx
	+ \alpha\int_{\frac{nt}{2}}^{nt}\left|\frac{3x - 2nt - 1}{n-1}\right|dx + O(1)                                                                           \\
	                            & = \left(\frac{5\alpha t^2}{12} - t + \frac{1}{\alpha}\right)(n - 1) + O(1),
\end{align*}
noting that we here have used $m = nt$, $t \in [0,1]$ and $\alpha$ constant with $t < \frac{1}{\alpha}$.
We conclude that for $n$ large enough, for $t = 0.216$ and for $\alpha = 3.29$ we have:
\begin{align}
	 & \sum_{i=1}^{n} \left( w(R_i) + w(C_i) \right) + \sum_{k=1}^{2n-1} \left( w(D_{k}^{-}) + w(D_{k}^{+}) \right) \label{eq:w_sum_final}                \\
	 & = \left(\frac{\alpha}{4}(1-t)^2 + \frac{1}{\alpha} - t + \alpha \frac{t^2}{4} + \frac{5\alpha^2t^2 - 12\alpha t + 12}{12\alpha}\right)(n-1) + O(1) \\ 
	 & \stackrel{\textrm{n large enough}}{\leq} 0.78379(n-1) < 0.784n = (1-t)n = n - |Q^\prime|,
\end{align}
which by Proposition~\ref{prop:weighting_method} implies that $Q^\prime$ is not completable.

For large enough $n$, we can always choose $m$ such that $m \; \equiv \; 1 \; (\textrm{mod} \; 6)$ and $t = 0.216 + \epsilon$ for $\epsilon \leq \frac{5}{n}$,
since $t = \frac{m}{n}$ and there are at most $5$ additional queens needed to have $m \equiv 1 \; (\textrm{mod} \; 6)$. The function
\begin{align*}
	f(t | \alpha) = \frac{\alpha}{4}(1-t)^2 + \frac{1}{\alpha} - t + \alpha \frac{t^2}{4} + \frac{5\alpha^2t^2 - 12\alpha t + 12}{12\alpha}
\end{align*}
being continuous in $t$ for all valid values of the parameter $\alpha$ implies that having $\epsilon$ small enough, equivalently $n$ large enough,
would leave the inequalities in (\ref{eq:w_sum_final}) unchanged, thus proving that this construction is not completable either.

Finally, for any odd $n$, embed the construction $Q^\prime$ into the upper-left $n \times n$ part of the
$(n+1) \times (n+1)$ board and give lines $R_{n+1}$ and $C_{n+1}$ weight $1$ to get the desired weighting
for even sized $n$.

Therefore, we improve the previously best known upper bound from \cite{Glock2022NQueens}.
\begin{theorem}
	For $n$ large enough $qc(n) < 0.216n$. \label{thm:improved_upperbound_qc}
\end{theorem}
Numerical experiments \cite{nielsen2025qclp} suggest that $t = 0.216$ is near optimal for
this $Q^\prime$ construction under the line-weighting method of \cite{Glock2022NQueens}.
However, the method may admit refinements, e.g., aggregating only the $n - m$
largest weights per line type, or incorporating additional row/column constraints, to further reduce the total weight.
It would be interesting to explore these methods to further improve upon the upper bound in Theorem~\ref{thm:improved_upperbound_qc}.

\backmatter

\bmhead{Acknowledgments}
I am grateful to Jakob Lemvig and Christian Henriksen for their supervision,
feedback on methodology, and helpful comments.

\bibliography{sn-bibliography}

\end{document}